\documentclass[12pt,preprintnumbers,amsmath,amssymb]{revtex4}
\usepackage{graphicx}% Include figure files
\usepackage{dcolumn}% Align table columns on decimal point
\usepackage{bm}% bold math
\usepackage{color}

\begin{document}

\title{Infinity increasing steps on the Collatz sequences}

\author{Eduardo M. K. Souza}
\affiliation{Departamento de Fisica, Universidade Federal de Sergipe 49100-000, Sao Cristovao SE, Brazil}

\date{\today}

\begin{abstract}
We intend to contribute to the Collatz dynamics problem by seeking to analyze the Collatz conjecture from the tree of numbers sequences. First, we show numerically that the distribution of odd numbers has an initial transient, and proceeds to a power law growth to its maximum. Second, using the formulation that uses only odd numbers, we present analytically a set of odd number sequences that is always increasing and can have an infinite number of terms.
\end{abstract}

\maketitle

\section{Introduction}

The Collatz's conjecture (CC) is related to a simple dynamics problem of the number theory in which starting from a positive natural number it always arrives at a final situation associated with a unique periodic sequence of numbers \cite{lag12,andr,wir}. 
Its proof has challenged a large number of scientists for many years. 
Despite it is apparently very simple, its proof has not been established until the present day. 
Some recent advances have raised hope that a proof can be presented in a shorter period of time. 
In short, to prove that the CC is false it is necessary to find (1) a non-trivial cycle, different from the final situation presented by Collatz and/or (2) a divergent orbit (infinity). 

In this work, it is presented a set of orbits on the Collatz sequence, in which the sequence of numbers is increasing and can have an infinite number of terms,
before reaching the final convergence. It is not discussed the steps sequence convergence after the sequence presented here, which we believe is the final situation presented by Collatz.
 
Let us present the CC. Consider a variable $X$ that takes values from the set of natural numbers $\mathbb{N}^{*}=\{ 1, 2, 3,...\}$. This variable evolves from a simple rule, so that one can think that $X$ is subjected to a dynamic in which $X(t+1)=C(X(t))$, for $t$ also belongs to the set of natural numbers $\mathbb{N}^{*}$.
The Collatz dynamics $C(X)$ is expressed through the rule \cite{lag12}

\begin{equation} \label{eq1}
C(X) = \left\{
\begin{array}{cc}
3X+1 & \text{if $X$ is odd},\\
X/2 & \text{if $X$ is even}.
\end{array}
\right.
\end{equation}

\noindent
Therefore, every natural number $X(t=1)$ generates a sequence $X(1)$, $X(2)=C(X(1))$, $X(3)=C(X(2))=C(C(X(1)))=C^{2}(X(1))$, ...
For example, it is easy to see that for $X(1)=1$ generates the sequence $1,2,4,1,2,4,1,...$, which is known as the Collatz periodic sequence (CPS). We see that $X(1)=2$ generates the sequence $2,1,4,2,1,4,...$; $X(1)=3$ generates the sequence $3,10,5,16,8,4,2,1,4,2,1,...$ and $X(1)=4$ generates the sequence $4,2,1,4,2,1,...$. 
Every sequences above fall into the CPS.
The Collatz's conjectured affirms that independent of an initial value $X(1)$, after a time $t$ the dynamics will lead to $X(t)$ for the CPS. 
This is, for every $X(1) \in \mathbb{N}^{*}$,
it is true the sequence $X(1)$, $X(2)$,..., $4$, $2$, $1$, $4$, $2$, $1$, ... 

Despite not having had an exact solution for the CC, some partial results, mainly from numerical tools, have shown great progress towards obtaining the final solution.
Numerical results have demonstrated the validity of the conjecture for values from $X(1)$ up to $10^{20}$ (2017). It has been demonstrated that the existence of a cycle different from the Collatz one in which the sequence can arrive, must have a minimum length of 17,087,915 steps \cite{eli93}.

The rest of this paper is organized as it follows: Sec. II presents a numerical analysis about the distribution of odd numbers on Collatz dynamics. Section III discusses analytically the Collatz dynamics of odd numbers.
Sec. IV closes the paper with concluding remarks.

\section{Density of odd numbers}

From Eq. (\ref{eq1}), it is trivial to observe that $X(t+1) < X(t)$ for even $X(t)$ and $X(t+1) > X(t)$ for odd $X(t)$. 
Therefore, an interesting study of dynamics is to observe the even and odd sequence in dynamics of $X$.
As odd numbers lead to the sequence of larger numbers by being multiplied by a factor 3, while even numbers are divided by a factor 2, we should expect that the number of odd numbers in the sequences to be less than the number of even numbers, as we show below.

It is possible to define a parity sequence
\begin{equation}
\{ X(1),X(2),X(3),...\} \;\;\;\;\; \rightarrow \;\;\;\;\; \{ y(1),y(2),y(3),...\} ,
\end{equation} 
where $y(t)=0$ if $X(t)$ is even, and  $y(t)=1$ if $X(t)$ is odd \cite{eve}.
It is important to note that this problem is not the same as studying CC in base 2 \cite{ba2}.
In the present situation we have a sequence of even and odd numbers represented by zeros and ones. This does not mean that the numbers are represented in base 2.
We are interested in evaluating the trajectories of the quantities of even and odd numbers in the sequence \cite{eve,ter}. We do a numerical analysis of this problem. We define 
\begin{equation}
P_{odd} (X) \equiv \frac{ \sum_{t=1}^{\sigma_{\infty}(X)} y(t)}{ \sigma_{\infty}(X) },
\end{equation}
which gives the fraction of odd numbers regards the total steps need to $X$ goes to $1$, given by $\sigma_{\infty} (X) = \text{inf} \{ k: C^{k} (X)=1 \}$.

\begin{figure}
	\centering
	\includegraphics[width=10cm,angle=0]{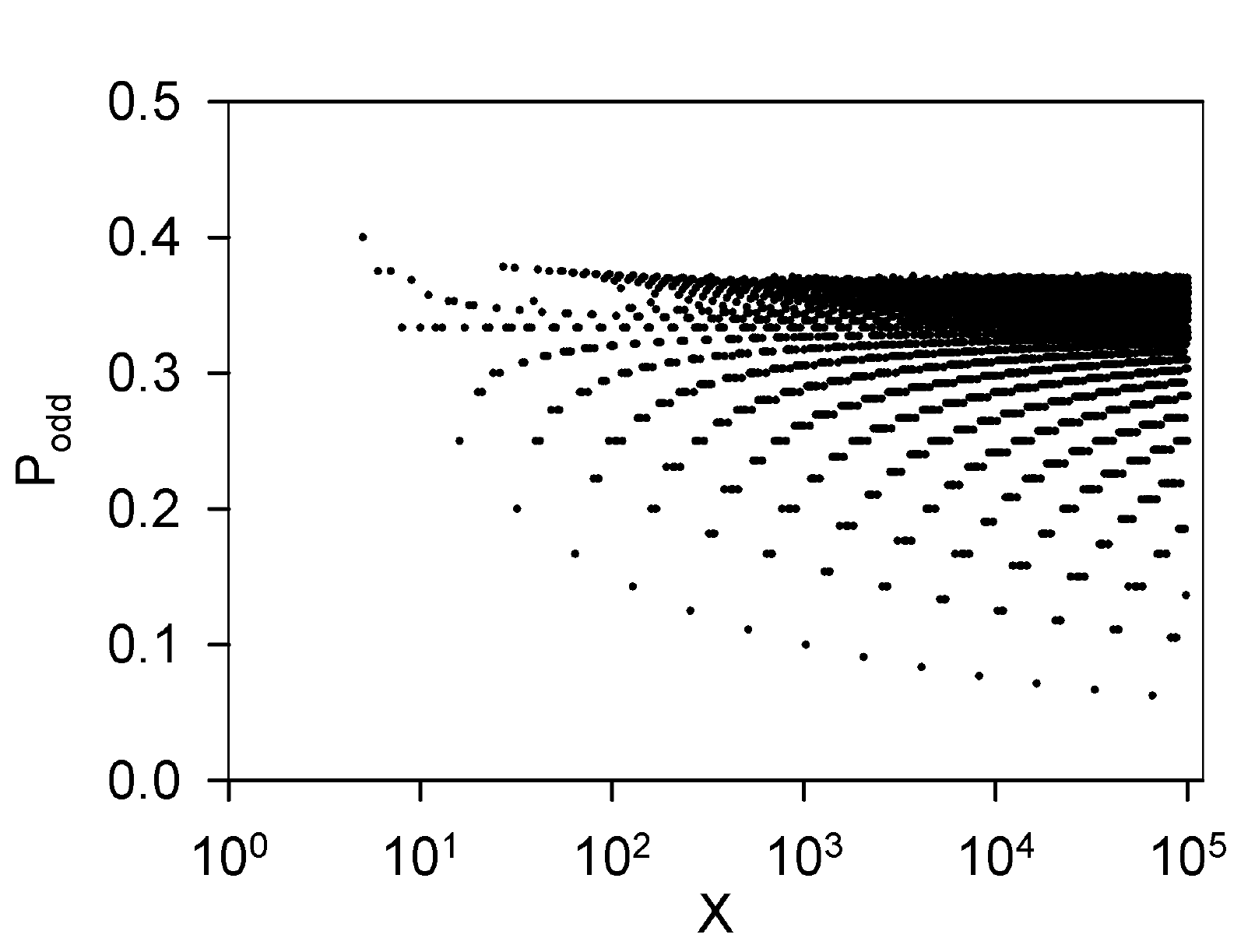}
	\caption{Fraction of odd numbers $P_{odd}$ regarding the total steps need to $X$ goes to $1$ as a function of $X$ up to $X \le 10^{5}$. } \label{poddx}
\end{figure}

Fig. \ref{poddx} shows $P_{odd}$ as a function of $X$ up to $10^{5}$. 
It is seen that for $X > 84$ we find $P_{odd}^{max} \approx 0.372$. 
The convergence of $P_{odd}^{max}$ is given in such a way that, for $X > 10^{2}$ we find $P_{odd} (X) \le 0.37190$, for $X > 10^{3}$ we find $P_{odd} (X) \le 0.37168$ and for $X > 10^{4}$ we find $P_{odd} (X) \le 0.37168$. 

\begin{figure}
	\centering
	\includegraphics[width=10cm,angle=0]{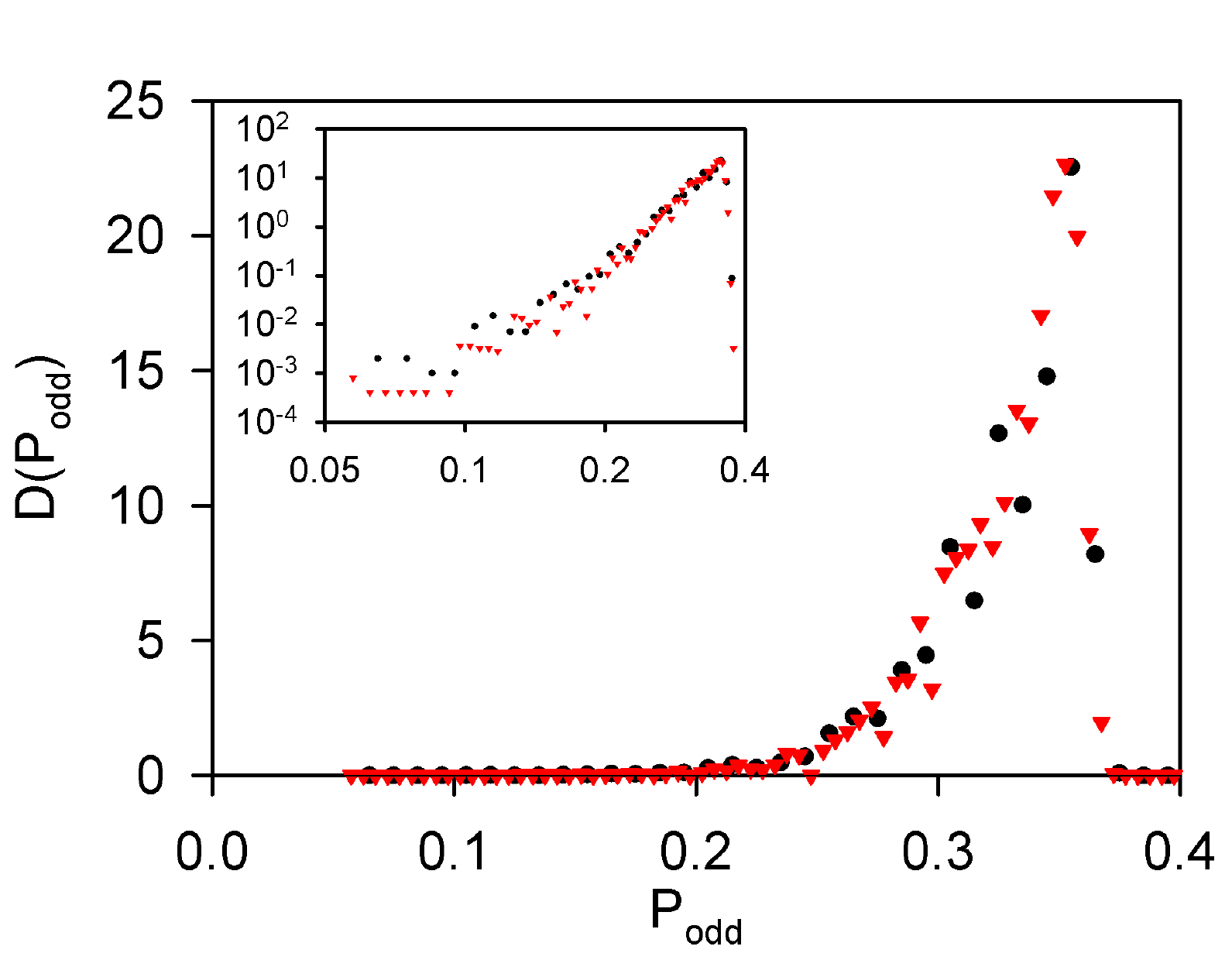}
	\caption{Probability distribution $D(P_{odd})$ of $P_{odd}$ for $X \le 10^{5}$ (black circle) and $X \le 5.10^{5}$ (red triangle). INSET: Same graph in $\log \times \log$ scale.} \label{dpodd}
\end{figure}

We can see in Fig. \ref{poddx} a dispersion of $P_{odd}$ values. Therefore, it is interesting to see the histogram of these $P_{odd}$. 
In Fig. \ref{dpodd} we can see the probability distribution of the $P_{odd}$ using the intervals $X \le 10^{5}$ (black circles) and $X \le 5.10^{5}$ (red triangles). 
In the inset we can see the same graph in $\log \times \log$ scale.
 Interesting to see that the distribution of odd numbers $D=(P_{odd})$ has an initial transient, and proceeds to a power law growth ($D=(P_{odd})^{\alpha}$ with $\alpha \sim 9$) to its maximum  $P_{odd}^{max}$ and then drops abruptly to zero.

From Fig. \ref{poddx}, the increase of $X$ brings a small increase in the growth slope, indicating a greater relevance for the $P_{odd}^{max}$ with the increase of $X$.
Our results are in agreement with previous ones in which $P_{odd} < 0.5$, therefore, the odd number is always less than the even number in the Collatz sequence. Our main result of this section is to have shown that from the pattern of Fig. \ref{poddx} the distribution of $P_{odd}$ follows a power law growth as shown in Fig. \ref{dpodd}.

\section{Collatz dynamics of odd numbers}

\begin{figure}
\centering
\includegraphics[width=12cm,angle=0]{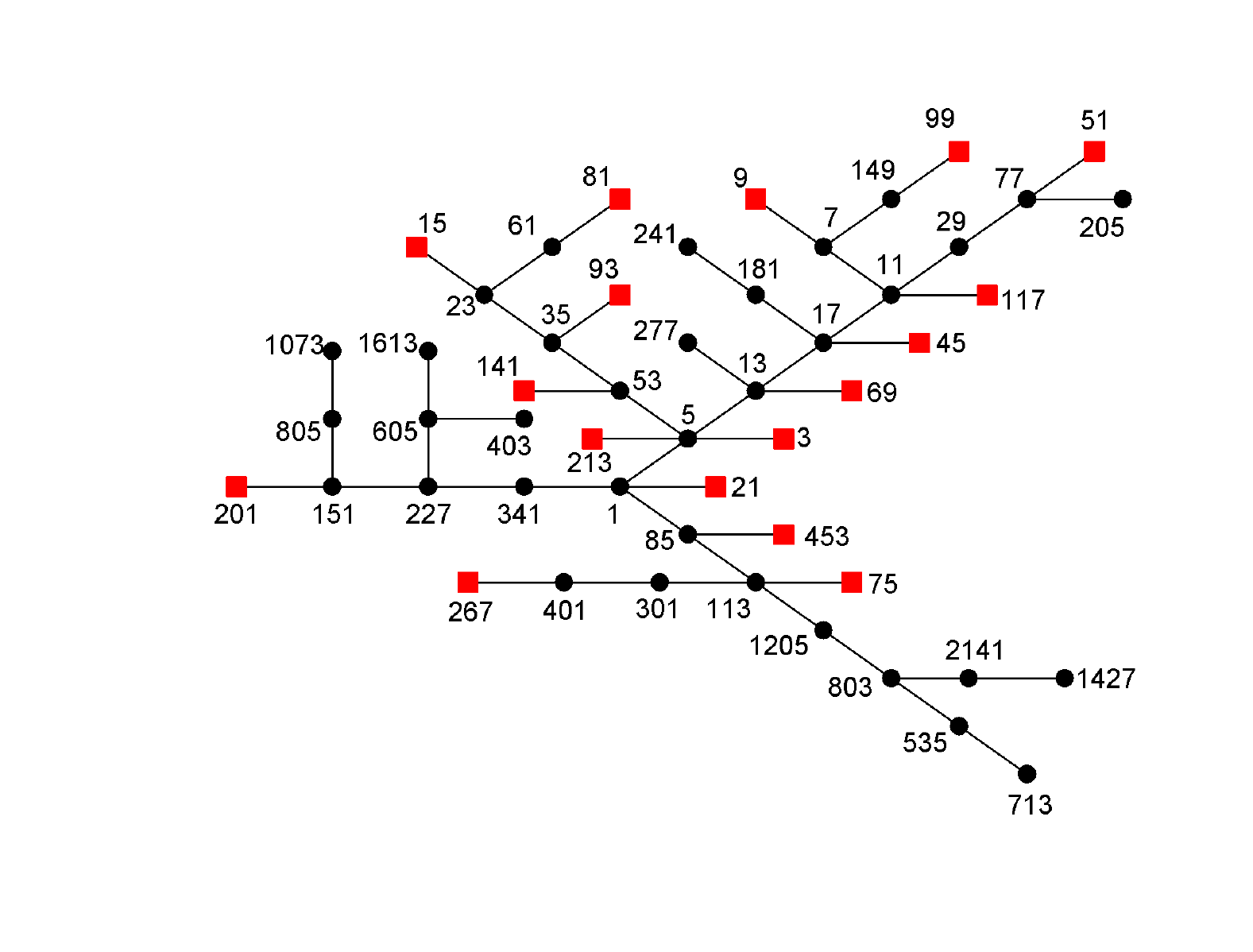}
\caption{Trees evolution of the Collatz dynamics showing only the odd numbers. 
Branches always start from a number multiple of $3$, represented by red squares.
Branches that start with black dots indicate that they are incomplete.} 
\label{tree}
\end{figure}

An alternative form of presenting Collatz dynamics that removes every even number is written as
$X_{odd}(t+1)=F(X_{odd}(t))$, where \cite{cham}
\begin{equation} 
	F(X_{odd})= \frac{3X_{odd}+1}{2^{m(3X+1)}} \;\;\;\;\;\; F: \mathbb{N}_{odd} \rightarrow \mathbb{N}_{odd}	
\label{cof3}
\end{equation}
which $m(X)$ is the number of factor $2$ contained in $X$. 
In this case we have a sequence of only odd numbers.
Fig. \ref{tree} illustrates the trees connection of odd numbers. This figure can be easily obtained from dynamics obtained by Eq. (\ref{cof3}).
The tree is formed by a sequence of branches that always lead to the number $1$. 
The branches always start from a number multiple of $3$ (represented by red squares).
There are no odd numbers in the sequence before odd numbers that are multiples of $3$ (see proof in Appendix).
The branches that start with black dots in Fig. \ref{tree} indicate that they are incomplete, as they are so large that they are not represented in the figure.

We can separate the trees of the dynamics of odd numbers into two sets.
The first set, we consider the cases that $X_{odd}(t+1)$ is always less than $X_{odd}(t)$.
The second set, we consider on the contrary, i.e., cases that $X_{odd}(t+1)$ is always greater than $X_{odd}(t)$.

The first set, cases that $X_{odd}(t+1) < X_{odd}(t)$, is associated with the odd numbers of the form $X_{odd}(t)=4n+1=5,9,13,17,...$, where
$n \in \mathbb{N}^{*}$. 
After some straightforward calculations, we found two subsets of these dynamics: 
i) the $X_{odd}(t)=(2^{p})n - [(1+2^{(p-1)})/3]$  goes to $ X_{odd}(t+1)=6n-1$ after $p=4,6,8,10,..$ steps; 
ii) the $X_{odd}(t)=(2^{p})n - [(1+5.2^{(p-1)})/3]$ goes to $X_{odd}(t+1)=6n-5$ after $p=3,5,7,9,..$ steps. 
Therefore, between $X_{odd}(t)$ and $X_{odd}(t+1)$ the number of 
$X_{even}$ is $p-1$, therefore, the factor $2$ contained 
between $X_{odd}(t)$ and $X_{odd}(t+1)$ is $m(X)=p-1$.
Thus, we can write that, if
\begin{equation} \label{dodd}
X_{odd}(t) = \left\{
\begin{array}{ll}
(2^{p})n - [(1+2^{(p-1)})/3] & \text{for $p \in \mathbb{N}_{even} \ge 4$},\\
(2^{p})n - [(1+5 \times 2^{(p-1)})/3] & \text{for $p \in \mathbb{N}_{odd} \ge 3$}.
\end{array}
\right.
\end{equation}
than
\begin{equation} \label{dodd2}
X_{odd}(t+1) = \left\{
\begin{array}{c}
6n-1 \\
6n-5 
\end{array}
\right.
\end{equation}
where $n \in \mathbb{N}^{*}$.
For example, we can find the odd numbers greater than $5$ that reach the number $5$. 
From Eq. (\ref{dodd2}), we need to use $X_{odd}(t+1)=6n-1=5$, thus $n=1$ and $p$ must be even. Now, from Eq. (\ref{dodd}), we find $X_{odd}(t)=(2^{p}) - [(1+2^{(p-1}))/3]=13,53,213,853,...$ Fig. \ref{tree} shows the first three numbers (13, 53, 213) of the infinite possible dynamics that lead to the number five. 
The above example for the number five can be repeated for every other odd number,
which can be expressed by Eq. (\ref{dodd2}).

The second set, cases that $X_{odd}(t+1) > X_{odd}(t)$, is associated 
with the odd numbers $X_{odd}(t)=4n-1=3,7,11,15,...$, that always leads to the odd number of value, respectively, $X_{odd}(t+1)=5,11,17,23,...$.
Considering the sequence of even and odd numbers, we observed that two steps of the original Collatz sequence are always required for $X_{odd}(t)$ to reach $X_{odd}(t+1)$.
As there are two steps, it is easy to see that the dynamics is given by $X_{odd}(t+1)=C(C(X_{odd}(t)))$, first for an odd number and then for an even number, so that, if $X_{odd}=4n-1$ then $X_{odd}(t+1)=(3X_{odd}(t)+1)/2=(3(4n-1)+1)/2=6n-1$.
Hence, we can write
\begin{equation} \label{eob}
\text{if} \; X_{odd}(t)=4n-1 \; \text{than} \; X_{odd}(t+1)=6n-1 \; \text{for} \; n \in \mathbb{N}^{*}.  
\end{equation}

From the above rules, Eqs. (\ref{dodd}-\ref{eob}), it is possible to obtain directly all the sequences of odd numbers. 
From Eqs. (\ref{dodd}-\ref{eob}), it is straightforward to find that $X_{odd}(t)$ when it obeys the form  $X_{odd}(t)=2^{s+2}n-(2^{s+1}+1)$ generates a sequence such that
$X_{odd}(t) < X_{odd}(t+1) < X_{odd}(t+2) < ... < X_{odd}(t+s)$, where 
$X_{odd}(t+s)=F^{(s)}(X_{odd}(t))$.
This means a dynamic in which the sequence of $X_{odd}$ is increasing by $s$ steps from odd to odd numbers, which $s \; \text{and} \; n \in \mathbb{N}^{*}$.

From a slightly more elaborate calculation, we can find, in addition to the number of steps, the entire sequence of $X_{odd}(t)$ that increases in value.
Let us consider the variable $q \in \mathbb{N}$. We set up the Table \ref{tab} with the relationship between $q$ and $n$,

\begin{table}[t!]
	\caption{Relationship between $q$ and $n$}
	\label{tab}
\begin{center}
	\begin{tabular}{|c|ccccc|}
		q & 0 & 1 & 2 & 3 & 4 $\;$ ...\\ 
		\hline 
		& 1 & 2 & 4 & 8 & 16 ...\\ 
		&  & 3 & 6 & 12 & 24 ...\\ 
		n	&  &  & 9 & 18 & 36 ...\\ 
		&  &  &  & 27 & 54 ...\\ 
		&  &  &  &  & 81 ...\\
		&  &  &  &  & $\;\;\;$  ...\\  
	\end{tabular} 
\end{center}
\end{table}
Table \ref{tab} shows that, for a fixed $q$, we have $q+1$ terms of $n$. For example, if $q=2$ then $n=4$, $6$ and $9$. In the general case, it is easy to see that
the $q+1$ terms of $n$ are written as
\begin{equation} \label{seq}
n_{i} = 3^{0} 2^{q}, 3^{1} 2^{q-1}, 3^{2} 2^{q-2} ,..., 3^{i} 2^{q-i},..., 3^{q} 2^{0}, \;\; 0 < i < q.
\end{equation}
where $i$ symbolizes the different terms of $n$.  
Using the fact that we have an increasing sequence, it means that we are in the second set of the tree of dynamics of odd numbers. Therefore, using Eq. (\ref{eob}), we can find all the sequence of $X_{odd}(t)$, which is given by $X_{odd}(t_i)=4n_{i}-1$
for  $t_{i}=1,...i,...,q+1$. Finally, the final term has the form $X_{odd}(q+2)=6n_{q}-1$.
For example, from Eqs. (\ref{eob}) and (\ref{seq}) for $q=2$ we obtain $n_{i}=4,6,9$ and the three terms of the ascending sequence are $X_{odd}(1)=15$,  $X_{odd}(1)=23$ and $X_{odd}(3)=35$, and finally $X_{odd}(4)=53$.
We emphasize that $q$ is a generic integer, thus, it is possible to have an infinite Collatz sequence always increasing, considering that we can do $q \rightarrow \infty$.

We can further generalize these increasing sequences and consider all $p$ prime numbers except $2$ and $3$. 
In this way, for $p=1,5,7,11,13,...,$ we have
\begin{equation} \label{seqg}
X_{odd}(t_i)=4\tilde{n}_{i}-1 \;\;\; \text{for} \;\;\; \tilde{n}_{i}=p \times n_{i} 
\end{equation}
and $X_{odd}(q+2)=6\tilde{n}_{q}-1$.
This generalization is done only with prime numbers, since all other numbers are in the table of sequences of previous $n_{i}$.
Thus, we have infinite sequences that create $q+2$ terms of the Collatz sequence of increasing odd numbers, in which it is possible to consider $q \rightarrow \infty$.
This is the main result of this section.

\section{Conclusions}

To summarize, in this paper we present a numerical analysis to find that the fraction of odd numbers $P_{odd}$ regarding the total steps of the Collatz sequences. We obtain that the distribution of odd numbers $D=(P_{odd})$ has an initial transient, and proceeds to a power law growth to its maximum  $P_{odd}^{max} \approx 0.372$ and then drops abruptly to zero.

We also present a set of infinite sequences of odd numbers that are always increasing. We found the sequence of numbers and the number of terms. It is direct to see that it is possible to have infinite terms on the sequences. It is important to emphasize that we do not discuss the convergence \cite{davi} of the sequence of steps after the sequences presented in this work.

\section*{Acknowledgements}
I am greatly indebted with A.M.C. de Souza for extremely valuable discussions and remarks, and L. K. Souza for a careful reading of the manuscript.

\section*{Appendix}

From Eqs. (\ref{dodd2}) and (\ref{eob}), we find that for all $X_{odd}$ we must have
only $X_{odd}(t+1)=6n-1$ or $X_{odd}(t+1)=6n-5$. Therefore, $X_{odd}(t+1)=2 \; \text{mod(3)}$ or $X_{odd}(t+1)=1 \; \text{mod(3)}$.
Thus,  $X_{odd}(t+1) \ne 0 \; \text{mod(3)}$, and $X_{odd}(t+1)$ never is a multiple of $3$.
Therefore, odd multiples of 3 never come from another odd number. That is, they are always the beginning of the branches of the odd-sequence tree, as represented in Fig. \ref{tree}. 


\begin{thebibliography}{00}

\bibitem{lag12} J. C. Lagarias {\it The ultimate challenge : the 3X+1 problem},  (ed. American Mathematical Society. 2012).

\bibitem{andr} S¸ Andrei, M. Kudlek, and R. S¸ Niculescu,
{\it Some results on the Collatz problem}, Acta Informatica {\bf 37}, 145–160 (2000).

\bibitem{wir} G. J. Wirsching, {\it The Dynamical System Generated by the 3n+1 Function}, (Lecture Notes in Math. 1681, Springer, 1998).

\bibitem{eli93} S. Eliahou,  {\it The the 3X+1 problem: new lower bounds on nontrivial cycle lengths}, Discrete Math. \textbf{118}, 45-56 (1993).

\bibitem{eve} C. J. Everett, {\it Iteration of the number-theoretic function}, in 
The ultimate challenge : the 3X+1 problem, Jeffrey C. Lagarias. editor (American Mathematical Society. 2012) p. 225. Originally published in Los Alamos report LA-6449—MS (1976).

\bibitem{ba2} J. Capco. {\it Odd Collatz Sequence and Binary Representations},  https://hal.archives-ouvertes.fr/hal-02062503, hal-02062503 (2019).

\bibitem{ter} R. Terras, {\it  On the existence of a density}, Acta Arithimetica \textbf{35}, 101-102 (1979).

\bibitem{cham} M. Chamberland, {\it A 3x + 1 Survey: Number Theory and Dynamical Systems}, in The ultimate challenge : the 3X+1 problem, Jeffrey C. Lagarias. editor (American Mathematical Society. 2012) p. 57.

\bibitem{davi} D. Barina, {\it Convergence verification of the Collatz problem},
The Journal of Supercomputing {\bf 77}, 2681–2688 (2021).

\end{thebibliography}
\end{document}